\documentclass{amsart}
\usepackage{amssymb,euscript,amsmath, mathrsfs}
\usepackage{bm}
\usepackage{rotating}
\usepackage[table]{xcolor}

\makeindex

\newcounter{ENUM}

\newcommand{\margh}[1]{}

\input xy
\xyoption{all}
\CompileMatrices

\newcommand{\wt}[1]{\widetilde{#1}}

\def\ZZ{{\mathbb Z}}

\def\PP{{\mathbb P}}

\def\cM{{\mathcal M}}

\def\sL{{\mathscr L}}

\def\sO{{\mathscr O}}
\def\sV{{\mathscr V}}

\def\fg{{\mathfrak g}}

\def\Hom{\operatorname{Hom}}

\def\Spec{\operatorname{Spec}}
\def\Pic{\operatorname{Pic}}

\def\Gr{\operatorname{Gr}}

\def\rk{\operatorname{rk}}

\def\MR{\operatorname{MR}}
\def\codim{\operatorname{codim}}

\def\inv{\operatorname{inv}}
\def\Jac{\operatorname{Jac}}

\newtheorem{thm}{Theorem}[section]
\newtheorem{prop}[thm]{Proposition}

\newtheorem{cor}[thm]{Corollary}

\theoremstyle{definition}
\newtheorem{defn}[thm]{Definition}

\newtheorem{ex}[thm]{Example}
\newtheorem{sit}[thm]{Situation}
\theoremstyle{remark}
\newtheorem{notn}[thm]{Notation}
\newtheorem{rem}[thm]{Remark}

\numberwithin{equation}{section}

\begin{document}
\title{The Gieseker-Petri theorem and imposed ramification}
\author{Melody Chan}
\author{Brian Osserman}
\author{Nathan Pflueger}
\begin{abstract} We prove a smoothness result for spaces of linear series
with prescribed ramification on twice-marked elliptic curves. In
characteristic $0$, we then
apply the Eisenbud-Harris theory of limit linear series to deduce a new 
proof of the Gieseker-Petri theorem, along with a generalization to spaces 
of linear series with prescribed ramification at up to two points. Our
main calculation involves the intersection of two Schubert cycles in a
Grassmannian associated to almost-transverse flags.
\end{abstract}

\subjclass[2010]{14H51, 14M15}

\thanks{The first author is supported by the Henry Merritt Wriston 
Fellowship and by NSF DMS-1701924. The second author is partially supported 
by a grant from the Simons Foundation \#279151.}

\maketitle

\section{Introduction}
The classical Brill-Noether theorem states that if we are given 
$g,r,d \geq 0$, a general curve $X$ of genus $g$ carries a linear series
$(\sL,V)$ of projective dimension $r$ and degree $d$ if and only if the quantity
$$\rho(g,r,d):=g-(r+1)(r+g-d)$$
is nonnegative \cite{g-h1}. Moreover, in this case 
the moduli space $G^r_d(X)$ of such linear series has pure dimension
$\rho$. This statement was generalized by Eisenbud and Harris
to allow for imposed ramification: given marked points 
$P_1,\dots,P_n \in X$, and sequences $0\le a^i_0<\dots<a^i_r \leq d$ for
$i=1,\dots,n$, consider the moduli space 
$G^r_d(X,(P_1,a^1_{\bullet}),\dots,(P_n,a^n_{\bullet})) \subseteq
G^r_d(X)$ parametrizing linear series with vanishing sequence at least
$a^i_{\bullet}$ at each of the $P_i$. Then Eisenbud and Harris used their
theory of limit linear series to show in \cite{e-h1}
that in characteristic $0$, if $(X,P_1,\dots,P_n)$ is a general $n$-marked 
curve of genus $g$, the dimension of
$G^r_d(X,(P_1,a^1_{\bullet}),\dots,(P_n,a^n_{\bullet}))$---if it is 
nonempty---is given by the generalized formula
$$\rho(g,r,d,a^1_\bullet,\ldots,a^n_\bullet):=g-(r+1)(r+g-d)-\sum_{i=1}^n \sum_{j=0}^r (a^i_j-j).$$
The condition for nonemptiness is still combinatorial, but becomes more 
complicated in this context. 

This theorem fails in positive characteristic
for $n \geq 3$, but is still true if $n \leq 2$.  In this case, we also
have a simple criterion for nonemptiness. To state it, we shift
notation, supposing we have marked points $P,Q\in X$, and sequences
$a_{\bullet}$, $b_{\bullet}$. We then introduce the following notation:
$$\widehat{\rho}(g,r,d,a_\bullet,b_\bullet) :=
g-\sum_{j:a_j+b_{r-j} > d-g} a_j+b_{r-j}-(d-g).$$

We summarize what was previously known about the space 
$G^r_d(X,(P,a_{\bullet}),(Q,b_{\bullet}))$, as follows.

\begin{thm}\label{thm:bg}
Given $(g,r,d)$ nonnegative integers,
and sequences $0 \leq a_0<a_1 <\dots<a_r \leq d$,
$0 \leq b_0<b_1 <\dots<b_r \leq d$, let $(X,P,Q)$ be a twice-marked
smooth projective curve of genus $g$ over a field of any characteristic. 
Set $\rho=\rho(g,r,d,a_{\bullet},b_{\bullet})$ and
set $\widehat{\rho}=\widehat{\rho}(g,r,d,a_{\bullet},b_{\bullet})$. 

Suppose that $X$ and $P,Q$ are general. Then
$G^r_d(X,(P,a_{\bullet}),(Q,b_{\bullet}))$ is nonempty if and only if
$\widehat{\rho} \geq 0$, and if nonempty, has pure dimension
$\rho$. Furthermore, it is reduced and Cohen-Macaulay, and if 
$\widehat{\rho} \geq 1$, it is connected.
\end{thm}

For the nonemptiness and dimension statements, see \cite{os18};
for reducedness and connectedness, see \cite{os26}. The Cohen-Macaulayness
statement follows from the construction of 
$G^r_d(X,(P,a_{\bullet}),(Q,b_{\bullet}))$ (see for instance the proof
of Proposition \ref{prop:non-smooth} below) together with the 
Cohen-Macaulayness of relative Schubert cycles. 

What has remained open until now is the question of the singularities
of the space $G^r_d(X,(P,a_{\bullet}),(Q,b_{\bullet}))$. In the absence of marked points,
Gieseker in 1982 used degenerations to prove a conjecture of 
Petri that if $X$ is general, then the space $G^r_d(X)$ is also smooth 
\cite{gi1}. This proof was later simplified by Eisenbud and Harris
\cite{e-h2} and Welters \cite{we4} using ideas closely related to the
theory of limit linear series. These proofs all relied on proving injectivity 
of the
Petri map, by taking a hypothetical nonzero element of the kernel, and
carrying out a careful analysis of how it would behave under degeneration.
Recently, Jensen and Payne have given a different proof of the 
Gieseker-Petri theorem by tropical degeneration methods \cite{j-p1}.
We also mention that in Proposition 3.2 of \cite{c-h-t1}, Ciliberto,
Harris and Teixidor assert partial smoothness results in
the case of a single ramification point via an argument building on the
Eisenbud-Harris approach, but their statement omits a key hypothesis; 
see Remark \ref{rem:wrong} below. 

In this paper, we give a new proof of the Gieseker-Petri theorem, and
generalize it to the space $G^r_d(X,(P,a_{\bullet}),(Q,b_{\bullet}))$,
proving that the singular locus of this space consists precisely of linear series with a certain type of excess vanishing.  
Our Gieseker-Petri theorem with imposed ramification, Theorem~\ref{thm:main} below, generalizes two statements.
\begin{enumerate}
\item In the absence of marked points, it reduces to the Gieseker Petri theorem, which holds for curves of any genus.
\item With marked points allowed, but in genus $0$, it reduces to the well-known characterization of the singular loci of Schubert varieties and Richardson varieties.
\end{enumerate}

Indeed, in the case $g=0$, a single ramification condition corresponds to a
Schubert cycle in the Grassmannian $\mathrm{Gr}(r+1,\mathcal{O}_{\mathbb{P}^1}\!(d))$, while a pair of ramification 
conditions similarly corresponds to a Richardson variety. These spaces are
singular, and their singularities can be characterized precisely as
loci with a specific type of excess vanishing. Our main theorem
extends this characterization to all genera, and also deduces additional
consequences on the geometry of $G^r_d(X,(P,a_{\bullet}),(Q,b_{\bullet}))$.
To state it, the following preliminary notation will be helpful.

\begin{notn} If $(\sL,V)$ is a $\fg^r_d$ on a smooth projective curve
$X$, and $D$ is an effective divisor on $X$, write
$$V(-D)=V \cap \Gamma(X,\sL(-D)) \subseteq \Gamma(X,\sL).$$
\end{notn}

Thus, to say that $(\sL,V)$ has vanishing sequence at least $a_{\bullet}$  
at $P$ is equivalent to saying that 
\begin{equation}\label{eq:van-ineq} 
\dim V(-a_j P) \geq r+1-j
\end{equation}
for $j=0,\dots,r$.

We then make the following definition:

\begin{defn}\label{defn:good-open} In the situation of Theorem
\ref{thm:bg}, let
$$G^{r,\circ}_d(X,(P,a_{\bullet}),(Q,b_{\bullet})) \subseteq G^r_d(X,(P,a_{\bullet}),(Q,b_{\bullet}))$$
be the open subset consisting of $(\sL,V)$ such that
\eqref{eq:van-ineq} holds with equality for all $j>0$ such that 
$a_j>a_{j-1}+1$, and the analogous condition holds for  $(Q,b_\bullet)$.
\end{defn}

We see that $G^{r,\circ}_d(X,(P,a_{\bullet}),(Q,b_{\bullet}))$ 
contains all linear
series with precisely the prescribed vanishing at $P$ and $Q$, but it
also contains many linear series with more than the prescribed vanishing.
For instance, if $a_{\bullet} = b_{\bullet} = (0,1,\ldots,r)$ are both minimal, so that
$G^r_d(X,(P,a_{\bullet}),(Q,b_{\bullet}))=G^r_d(X)$, then we also have
$G^{r,\circ}_d(X,(P,a_{\bullet}),(Q,b_{\bullet}))=G^r_d(X)$.
Our main theorem is then the following.

\begin{thm}\label{thm:main} In the situation of Theorem \ref{thm:bg}, 
suppose further that we are in characteristic $0$. Then the smooth locus of
$G^r_d(X,(P,a_{\bullet}),(Q,b_{\bullet}))$ is precisely equal to 
$G^{r,\circ}_d(X,(P,a_{\bullet}),(Q,b_{\bullet}))$. 

Furthermore, the
space $G^r_d(X,(P,a_{\bullet}),(Q,b_{\bullet}))$ has singularities in 
codimension at least $3$, is normal, and when $\widehat{\rho} \geq 1$
is irreducible.
\end{thm}

Thus, we are in particular giving a new proof of the Gieseker-Petri
theorem (in characteristic $0$).
As an immediate consequence of Theorem~\ref{thm:main}, the twice-pointed Brill-Noether curves studied in \cite{c-l-p-t1}, as well as twice-pointed Brill-Noether surfaces \cite{c-p3, a-c-t1} are smooth.
The singular locus of $G^r_d(X,(P,a_{\bullet}),(Q,b_{\bullet}))$ can also be described using hook removals in Young diagrams; see Remark~\ref{rem:hook}.

Our proof proceeds by degenerating to a chain of elliptic curves,
and studying the geometry of the corresponding moduli space of
Eisenbud-Harris limit linear series. 
The key
idea in this step is that although the space of limit linear series
will be singular in codimension $1$, after base change and blowup one can 
ensure that any given point of 
$G^{r,\circ}_d(X,(P,a_{\bullet}),(Q,b_{\bullet}))$ on the generic fiber will
specialize to a smooth point of the limit linear series space of a chain of curves of genus 0 or 1. This
is where the characteristic-$0$ hypothesis comes in. The case of genus 0 is well-known, so our main 
calculation is  the following result, which does not depend on 
characteristic, concerning the case $g=1$.

\begin{thm}\label{thm:genus-1} In the situation of Theorem \ref{thm:bg}, 
suppose that $g=1$, and make the generality condition explicit as follows:
$X$ is arbitrary, and $P,Q$ are such that $P-Q$ is not a torsion
point of $\Pic^0(X)$ of order less than or equal to $d$.
Then the space $G^{r,\circ}_d(X,(P,a_{\bullet}),(Q,b_{\bullet}))$ is smooth.
\end{thm}

The proof of Theorem \ref{thm:genus-1} proceeds by consideration of
the morphism
\begin{equation}\label{eq:to-pic}
G^{r,\circ}_d(X,(P,a_{\bullet}),(Q,b_{\bullet})) \to \Pic^d(X).
\end{equation}
The main subtlety that needs to be addressed is that the map~\eqref{eq:to-pic} 
is not smooth. The fibers are each described as an intersection of a pair
of Schubert cycles in a Grassmannian.  But in finitely many fibers, namely the ones above line bundles of the form $\mathcal{O}_X(aP+(d\!-\!a)Q)$ for $0<a<d$, the
pairs of flags defining the Schubert cycles are not transverse, but only almost-transverse (see 
Definition~\ref{def:a-t}). We prove Theorem
\ref{thm:genus-1} by first showing that in fibers, the tangent spaces
have dimension at most $1$ greater than expected, and then showing
that at the points where the tangent space dimension jumps in the fiber,
there cannot be any horizontal tangent vectors. 

The statement on tangent spaces in fibers,
which is Corollary \ref{cor:at} below, takes place entirely inside the 
Grassmannian, and may be of independent interest. Indeed, \S 
\ref{sec:a-t} is a study of tangent spaces of intersections of
pairs of Schubert cycles, and we address the case of arbitrary pairs of
flags in Theorem \ref{thm:workhorse} and Remark \ref{rem:arbitrary-flags}.
We then prove Theorems \ref{thm:main} and \ref{thm:genus-1} in \S
\ref{sec:proofs}.

\section{Almost-transverse intersections of Schubert cycles}\label{sec:a-t}

It is well known that the intersection of two Schubert varieties associated
to transverse flags---commonly called a Richardson variety---is 
smooth on the open subset of points which are smooth in both Schubert 
varieties.  In this section we consider intersections of pairs of Schubert 
varieties associated to not necessarily transverse flags. Our analysis 
recovers the usual smoothness statement in the transverse case, but our main 
purpose is to analyze the almost-transverse case in Corollary \ref{cor:at}, 
where we characterize the smooth
points and show that the dimension of the tangent space jumps only by $1$ 
at the non-smooth points. While Schubert intersections and non-transverse 
flags have been studied by Vakil \cite{va7} and Coskun 
\cite{co3}, those situations involved studying the flat limits of transverse
intersections, rather than the direct analysis of the non-transverse
intersections required in the present work.\footnote{More precisely, they
work with closures of loci with prescribed behavior with respect to both
flags; these are in particular irreducible, so are not the same thing as
the intersection of Schubert cycles which we consider.}

We fix $k$ to be an algebraically closed field of any characteristic. 
Throughout this section, we will work entirely with $k$-valued (equivalently,
closed) points. 
We index our complete flags by codimension, so that for a complete flag $P^\bullet$ in a $d$-dimensional vector space $H$,
\begin{align*}
0 &= P^d \subset \cdots \subset P^1 \subset P^0 = H.
\end{align*} We fix further notation as
follows.

\begin{defn}\label{defn:schubert}
Given a $k$-vector space $H$ of finite dimension $d$ and 
a complete flag $P^{\bullet}$ in $H$, if we are given also
$a_{\bullet} = (a_0,\ldots,a_r)\in \ZZ^{r+1}$ 
with $$0 \le a_0 < \cdots < a_r < d,$$
we let $\Sigma_{P^{\bullet}, a_{\bullet}}$ be the \textbf{Schubert variety} 
defined as the closed subscheme of $\Gr(r+1,H)$ given by the set of 
$\Lambda \in \Gr(r+1,H)$ such that
\begin{equation}\label{eq:one-flag}
\dim (\Lambda\cap P^{a_i}) 
\ge r+1-i
\end{equation}
for $i=0,\ldots,r$.
\end{defn}

More precisely, the conditions in~\eqref{eq:one-flag} are determinantal, 
yielding a scheme structure on $\Sigma_{P^{\bullet},a_{\bullet}}$ (which turns out
to be reduced). In our notation, the codimension of $\Sigma_{P^{\bullet},a_{\bullet}}$
is given by $\sum_{i=0}^r (a_i\!-\!i)$.

\begin{defn}\label{defn:active}
With $a_{\bullet} = (a_0,\ldots , a_r)$ an increasing sequence as above, say that an index $i$ with $0\le i\le r$ is \textbf{active} in $a_{\bullet}$ if $i>0$ and 
$a_i > a_{i-1} + 1$, or $i=0$ and $a_0 >0.$
\end{defn}

\begin{defn}\label{defn:sigma-circ}
Let ${\Sigma}^\circ_{P^{\bullet},a_{\bullet}}$ be the open 
subscheme of 
$\Sigma_{P^{\bullet},a_{\bullet}}$ consisting of subspaces $\Lambda$ for which for every active index $i$, the 
inequality in~\eqref{eq:one-flag} is an equality. 
\end{defn}

\noindent Note that \eqref{eq:one-flag} is automatically an equality when $i=0$, so
in Definition \ref{defn:sigma-circ} we can restrict to positive active 
choices of $i$.

We fix the following situation throughout this section.

\begin{sit}\label{sit:flags} 
Let $H$ be a finite-dimensional $k$-vector space, and write $d:=\dim H$. Fix complete flags
$P^{\bullet}$, $Q^{\bullet}$ in $H$, and sequences
$a_{\bullet},b_{\bullet}\in \ZZ^{r+1}$ 
with $0 \le a_0 < \cdots < a_r < d$ and $0 \le b_0 < \cdots < b_r < d.$
\end{sit}

Recall that we are indexing by codimension; thus $\codim P^i=\codim Q^i=i$.   
Note that for any $\Lambda \in {\Sigma}_{P^{\bullet},a_{\bullet}}$, the distinct subspaces in the collection
$\Lambda\cap P^{j}$ form a complete flag in $\Lambda$; 
we denote the flag $\Lambda\cap P^{\bullet}$ by abuse of notation.  

We have the following description of the tangent space at any point in 
${\Sigma}_{P^{\bullet},a_{\bullet}}$. Tangent spaces to Schubert varieties are well 
understood \cite{b-l2}, but for the sake of completeness, we provide a 
description in the particular case that we need of Grassmannian Schubert 
varieties. 

\begin{prop}
\label{prop:tangent-schubert} 
\mbox{}
\begin{enumerate}
\item Given 
$\Lambda \in {\Sigma}_{P^{\bullet},a_{\bullet}}$, let $S$ be the set of active indices $i$ such that $\dim \Lambda\cap P^{a_i} = r+1-i$.  Then 
there is a canonical isomorphism of vector spaces
$$T_{\Lambda}\Sigma_{P^{\bullet},a_{\bullet}}\cong 
\left\{\phi\colon \Lambda\to H/\Lambda: \phi(\Lambda \cap P^{a_i}) 
\subseteq (P^{a_i}+\Lambda)/\Lambda \text{ for } i\in S\right\}.$$

\item In particular, the smooth locus of $\Sigma_{P^{\bullet},a_{\bullet}}$ is precisely
${\Sigma}^\circ_{P^{\bullet},a_{\bullet}}$.
\end{enumerate}
\end{prop}

\begin{proof}
By definition, $\Sigma_{P^{\bullet},a_{\bullet}}$ is the scheme-theoretic intersection 
of the following subschemes of $\Gr(r+1,H)$ (for $i=0,1,\cdots,r$):
$$\Sigma_i = \{ \Lambda \in \Gr(r+1,H)\colon \dim(\Lambda \cap P^{a_i}) \geq
r+1-i \}.$$
Define $\Sigma_i^\circ$ to be the open subscheme of $\Sigma_i$ where equality
holds. Note also that in fact
$\Sigma_{P^{\bullet},a_{\bullet}}$ can be cut out as the intersection of the
$\Sigma_i$ over all active indices $i$: this is immediate set-theoretically,
and is also true scheme-theoretically because whenever $a_{i+1}=a_i+1$,
the condition for $a_{i+1}$ is obtained from that of $a_i$ by adding a 
single row to the local matrix expression, and considering minors of
size one larger. Thus, every minor occuring in the $a_{i+1}$ condition can be 
expanded in terms of minors occuring in the $a_i$ condition.

The first statement of the proposition then follows immediately from the
following claim. For a fixed index $i$,
$$ T_{\Lambda} \Sigma_i =
\begin{cases}
\{\phi\colon \Lambda\to H/\Lambda: \phi(\Lambda \cap P^{a_i}) \subseteq
(P^{a_i}+\Lambda)/\Lambda \} & \mbox{ if $\Lambda \in \Sigma^\circ_i$}\\
T_{\Lambda} \Gr(r+1,H) & \mbox{ otherwise,}
\end{cases}$$
where we identify $T_\Lambda \Gr(r+1,H)$ with $\Hom(\Lambda,H/\Lambda)$ as
usual.

To prove this claim, one may work on an affine open subset of $\Gr(r+1,H)$, as
follows. Choose a basis of $H$ extending a basis of $\Lambda$; then an affine
neighborhood of $\Lambda$ is given by the set of $(r+1) \times d$ matrices
whose first $(r+1)$ columns form the identity matrix (where the point in
$\Gr(r+1,H)$ is given by taking the span of the rows). More precisely, for any
$k$-algebra $R$, we may identify the $R$-points of this open subscheme with
$R$-valued matrices whose first $r+1$ columns form the identity matrix. In
particular, taking $R = k[\epsilon]$, the tangent space $T_{\Lambda}
\Gr(r+1,H)$ is identified with matrices in block form $( I\  \epsilon M)$,
where $M$ is a matrix of values of $k$; the matrix $M$ then determines an
element of $\Hom(\Lambda,H / \Lambda)$. Now, we may further assume that the
chosen basis of $H$ also includes a basis of $P^{a_i}$ as a subset. Then the
$R$-points of $\Sigma_i$ consist of those matrices such that the submatrix
consisting of all columns not corresponding to the basis of $P^{a_i}$ has rank
at most $i$. Assuming that we order our basis of $\Lambda$ so that a basis of $\Lambda
\cap P^{a_i}$ comes at the end, the submatrix in
question has the form
$$\begin{pmatrix} I & A \\ 0 & B \end{pmatrix},$$
where the size of the identity matrix in the upper left is $\dim (\Lambda /
(\Lambda \cap P^{a_i}))$. Therefore, lying in $\Sigma_i$ corresponds to the
condition that
$$\rk(B) \leq i - (r+1) + \dim (\Lambda \cap P^{a_i}).$$
Now, specialize to the case $R = k[\epsilon]$, and consider a tangent vector to
$\Gr(r+1,H)$ at $\Lambda$. The submatrix $B$ is a multiple of $\epsilon$.
Therefore all $2 \times 2$ and larger minors of $B$ are guaranteed to vanish.
Thus in the case $\dim (\Lambda \cap P^{a_i}) > r+1-i$ (i.e. $\Lambda \not\in
\Sigma_i^\circ$), all tangent vectors to $\Gr(r+1,H)$ at $\Lambda$ are also
tangent vectors to $\Sigma_i$ at $\Lambda$. On the other hand, when $\Lambda
\in \Sigma_i^\circ$, a tangent vector to $\Gr(r+1,H)$ at $\Lambda$ is a tangent
vector to $\Sigma_i$ if and only if the matrix $B$ vanishes entirely. This
condition can be made intrinsic by observing that, if $\phi\colon \Lambda
\rightarrow H / \Lambda$ is the linear map encoding a tangent vector, then $B$
is a matrix representation for the linear map $\Lambda \cap P^{a_i} \to H /
(P^{a_i} + \Lambda)$ induced by $\phi$. Therefore it follows that, in the case
$\Lambda \in \Sigma_i^\circ$, $\phi$ described a tangent vector to $\Sigma_i$
if and only if $\phi(\Lambda \cap P^{a_i}) \subseteq (P^{a_i} +
\Lambda)/\Lambda$. This proves the claim, and the first statement of
the proposition.

The second statement follows by direct computation of the codimension
imposed by the conditions on the tangent space in the first part.
If we have $i \in S$, let $i_n$ denote the next (greater) element of $S$,
setting $i_n=r+1$ if $i$ is maximal in $S$. By starting from the condition
imposed at the maximal element of $S$, and inductively working downwards,
one computes that the codimension of the tangent space is given by
$$\sum_{i \in S} (i_n-i)(a_i-i).$$
Each term of this sum is always less than or equal to 
$\sum_{j=i}^{i_n-1} (a_j-j)$, with equality if and only if there are no
actives indices strictly between $i$ and $i_n$. The proposition follows.
\end{proof}

\begin{cor}
\label{cor:tangent-schubert-circ} Given 
$\Lambda \in {\Sigma}^\circ_{P^{\bullet},a_{\bullet}}$,
there is a canonical isomorphism of vector spaces
$$T_{\Lambda}\Sigma_{P^{\bullet},a_{\bullet}}\cong 
\left\{\phi\colon \Lambda\to H/\Lambda: \phi(\Lambda \cap P^{a_i}) 
\subseteq (P^{a_i}+\Lambda)/\Lambda \text{ for active } i=0,\ldots,r\right\}.$$
\end{cor}

%
%

Following Definition 4.1 of \cite{c-p3}, we define:

\begin{defn} \label{def:a-t} Two complete flags $P^{\bullet}$ and $Q^{\bullet}$ are called 
\textbf{almost-transverse} if there exists an index $t\in\{1,\ldots,d-1\}$ 
such that
$$\dim P^i \cap Q^{d-i} 
= \begin{cases} 0 & \text{ if $i\ne t$,}\\ 1 & \text{ if }i=t.\end{cases}$$
\end{defn}

More generally, we have the following statement, which is easy to check:

\begin{prop}\label{prop:sigma} 
There is a unique permutation $\sigma \in S_d$ associated to the
flags $P^{\bullet}$ and $Q^{\bullet}$ with the property 
that there exists a basis $e_1,\ldots,e_d$ for $H$ satisfying
$$e_i \in P^{i-1}\smallsetminus P^i, \quad \text{and} \quad
e_i \in Q^{\sigma(i)-1} \smallsetminus Q^{\sigma(i)}.$$

Such a basis can also be characterized by the property that for all indices 
$i$ and $j$,
$P^i\cap Q^j$ is spanned by $\{e_1,\ldots,e_d\}\cap P^i\cap Q^j$.
In particular, if $\dim P^i\cap Q^j = 1$ then $P^i\cap Q^j$ contains one
of the $e_{\ell}$.
\end{prop}

\begin{defn} We refer to a basis as in Proposition \ref{prop:sigma} as
a $(P^{\bullet},Q^{\bullet})$-\textbf{basis}. 
\end{defn}

Thus, following the notation of Proposition \ref{prop:sigma}, we have
that $P^{\bullet}=Q^{\bullet}$ if and only if $\sigma=\mathrm{id}$, 
$P^{\bullet}$ and $Q^{\bullet}$ are transverse if and only if 
$\sigma = \omega := (d,d-1,\ldots,1)$, and $P^{\bullet}$ and $Q^{\bullet}$ 
are almost-transverse if and only if $\sigma$ is the composition of $\omega$ 
with an adjacent transposition.

We set $\rho = \rho(1,r,d,a_{\bullet},b_{\bullet})$, so that 
$$\rho-1=(r+1)(d-r-1)-\sum_{i=0}^r (a_i-i)-\sum_{i=0}^r (b_i-i)
= \rho(0,r,d-1,a_{\bullet},b_{\bullet})$$
is precisely the expected dimension of 
$\Sigma_{P^{\bullet},a_{\bullet}} \cap \Sigma_{Q^{\bullet},b_{\bullet}}$.
Also recall the definition of the complete flags $\Lambda\cap P^\bullet$ and $\Lambda\cap Q^\bullet$ from Situation~\ref{sit:flags}.

\begin{thm}\label{thm:workhorse}
Given $\Lambda 
\in {\Sigma}^\circ_{P^{\bullet},a_{\bullet}} \cap {\Sigma}^\circ_{Q^{\bullet},b_{\bullet}}$,
let $\sigma \in S_{r+1}$ denote the permutation associated to
$\Lambda\cap P^{\bullet}$ and $\Lambda\cap Q^{\bullet}$ in $\Lambda$
by Proposition \ref{prop:sigma}. 
Given any $j\in \{0,\ldots,r\}$, let 
$$\mu(j) = \max\{i\le j: i\text{ is active in $a_{\bullet}$}\},$$
setting $\mu(j) = 0$ if no such $i$ exists.  Set $m(j) = a_{\mu(j)}.$
Similarly, let 
$$\nu(j) = \max\{i\le \sigma(j): i\text{ is active in $b_{\bullet}$}\},$$ 
setting $\nu(j) = 0$ if no such $j$ exists. Set $n(j) = b_{\nu(j)}$.
Then
\begin{equation}\label{eq:the-answer}
\dim \, 
T_\Lambda (\Sigma_{P^{\bullet},a_{\bullet}} \cap \Sigma_{Q^{\bullet},b_{\bullet}})
= \rho-1 + \sum_{j=0}^r \codim_H (P^{m(j)} + Q^{n(j)} + \Lambda).
\end{equation}
\end{thm}

\begin{proof}
Let $\lambda_0,\ldots,\lambda_r$ be a 
$(\Lambda\cap P^{\bullet},\Lambda\cap Q^{\bullet})$-basis for 
$\Lambda$.   Then for any $i$ active in $a_{\bullet}$, respectively $b_{\bullet}$, have
\begin{equation}\label{eq:lambda-basis}
\Lambda \cap P^{a_i} =\langle \lambda_i,\ldots,\lambda_r\rangle,\qquad \Lambda\cap Q^{b_i} = \langle \lambda_{\sigma^{-1}(i)}, \ldots,\lambda_{\sigma^{-1}(r)}\rangle.
\end{equation}
%
In other words, given any $j$, and any $i$ that is active in $a_{\bullet}$, we have $\lambda_j\in \Lambda\cap P^{a_i}$ if and only if $i\le j$.  Similarly, for any $i$ that is active in $b_{\bullet}$, we have $\lambda_j\in \Lambda\cap Q^{b_i}$ if and only if $i\le\sigma( j)$.  
By Corollary~\ref{cor:tangent-schubert-circ}, we have isomorphisms
\begin{eqnarray*}
T_\Lambda (\Sigma_{P^{\bullet},a_{\bullet}} \cap \Sigma_{Q^{\bullet},b_{\bullet}})
 &\cong& \{\phi\colon \Lambda\to H/\Lambda ~:~ \phi(\lambda_j) \in (P^{m(j)}+\Lambda)/\Lambda \cap (Q^{n(j)}+\Lambda)/\Lambda\}\\
 &\cong&  \bigoplus_{j=0}^r \mathrm{Hom} \left(\langle \lambda_j\rangle, (P^{m(j)}+\Lambda)/\Lambda \cap (Q^{n(j)}+\Lambda)/\Lambda\right).\end{eqnarray*}
 We are thus reduced to computing the dimensions $(P^{m(j)}+\Lambda)/\Lambda \cap (Q^{n(j)}+\Lambda)/\Lambda$, which are equal to
$$\dim (P^{m(j)}+\Lambda)+ \dim (Q^{n(j)}+\Lambda)
- \dim (P^{m(j)}+ Q^{n(j)}+\Lambda)-\dim \Lambda.$$
 Moreover, the first two terms are determined by the fact that $\dim P^{m(j)}\cap \Lambda = r+1-\mu(j)$ and $\dim Q^{n(j)} \cap \Lambda = r+1-\nu(j)$, by the assumption that $\mu(j)$ and $\nu(j)$ are active or are equal to $0$. Furthermore note by definition of $\mu$ and $\nu$ that $$a_j - j = m(j) - \mu(j),\qquad b_{\sigma(j)} - \sigma(j) = n(j) - \nu(j).$$
 A straightforward calculation then produces~\eqref{eq:the-answer}.
\end{proof}

We observe that the well-known case of transverse flags follows
immediately from Theorem \ref{thm:workhorse}.

\begin{cor}\label{cor:transverse}
If $P^{\bullet}$ and $Q^{\bullet}$ are transverse, then 
$\dim \, T_\Lambda (\Sigma_{P^{\bullet},a_{\bullet}} \cap \Sigma_{Q^{\bullet},b_{\bullet}})
= \rho-1 $ for all
$\Lambda 
\in {\Sigma}^\circ_{P^{\bullet},a_{\bullet}} \cap {\Sigma}^\circ_{Q^{\bullet},b_{\bullet}}$.
\end{cor}

\begin{proof} Following the notation of the proof of Theorem
\ref{thm:workhorse}, for each $j$ we have 
$\lambda_j  \in P^{m(j)} \cap Q^{n(j)} $ by construction. 
Since $P^{\bullet}$ and $Q^{\bullet}$ are transverse, it follows that 
$m(j) + n(j) < d$, so 
$P^{m(j)} + Q^{n(j)} + \Lambda = P^{m(j)} + Q^{n(j)} = H$.
\end{proof}

More importantly, we can also deduce the desired statement in the
almost-transverse case.

\begin{cor}\label{cor:at}
Given $\Lambda 
\in {\Sigma}^\circ_{P^{\bullet},a_{\bullet}} \cap {\Sigma}^\circ_{Q^{\bullet},b_{\bullet}}$,
suppose $P^{\bullet}$ and $Q^{\bullet}$ are almost-transverse, with $t+t'=d$ 
such that $\dim P^t \cap Q^{t'} =1$.  

Suppose first that $t=a_i$ for $i$ active in $a_{\bullet}$, that $t'=b_{i'}$ for $i'$ active in $b_{\bullet}$, and that 
$$P^t\cap Q^{t'} \subseteq \Lambda \subseteq P^t + Q^{t'}.$$ 
Then $$\dim \, 
T_\Lambda (\Sigma_{P^{\bullet},a_{\bullet}} \cap \Sigma_{Q^{\bullet},b_{\bullet}})
=\rho.$$
If those conditions do not all hold, then $$\dim \, 
T_\Lambda (\Sigma_{P^{\bullet},a_{\bullet}} \cap \Sigma_{Q^{\bullet},b_{\bullet}})
=\rho-1.$$
\end{cor}

\begin{proof}
Let $\lambda_0,\ldots,\lambda_r$ be a 
$(\Lambda\cap P^{\bullet},\Lambda\cap Q^{\bullet})$-basis for 
$\Lambda$.
First suppose that  $t=a_i$ for $i$ active in $a_{\bullet}$, and $t'=b_{i'}$ for $i'$ active in $b_{\bullet}$, and that $P^t\cap Q^{t'} \subseteq \Lambda \subseteq P^t + Q^{t'}.$  We will deduce that $\dim \, 
T_\Lambda (\Sigma_{P^{\bullet},a_{\bullet}} \cap \Sigma_{Q^{\bullet},b_{\bullet}})
=\rho.$   

We have that $\Lambda\cap P^t$ and $\Lambda\cap Q^{t'}$ are elements in the flags $\Lambda\cap P^\bullet$ and $\Lambda\cap Q^\bullet$ respectively with intersection $\Lambda\cap P^t \cap Q^{t'}$ of dimension 1.  Proposition~\ref{prop:sigma} implies that $P^t \cap Q^{t'} = \langle \lambda_j \rangle$ for a unique $j\in\{0,\ldots,r\}$.   By Theorem~\ref{thm:workhorse} it is enough to show that for each $j'\in\{0,\ldots,r\}$, $$\codim_H (P^{m(j')} + Q^{n(j')} + \Lambda) = \begin{cases}
1 &\text{ if }j'=j,\\
0&\text{ if }j'\ne j.
\end{cases}$$

Now, for $j'\ne j$, the fact that $\lambda_{j'}\in P^{m(j')}\cap Q^{n(j')}$ and $P^\bullet$ and $Q^\bullet$ are almost-transverse implies that either $m(j')+n(j')<d$, or that $m(j') = t$ and $n(j') = t'$.  But the latter case cannot be, since then both $\lambda_j, \lambda_{j'}\in P^{t}\cap Q^{t'}$, contradicting that $\dim P^{t}\cap Q^{t'}=1$.  Therefore $m(j')+n(j')<d$ and $$P^{m(j')}+ Q^{n(j')} = P^{m(j')}+ Q^{n(j')}+\Lambda = H,$$ as desired. 

Next, to show that $\codim_H (P^{m(j)} + Q^{n(j)}+\Lambda )=1$, we claim that $m(j) = t$ and $n(j) = t'$. Recall that $a_i = t$ and $a_{i'} =t'$. By assumption, $i$ is active in $a_{\bullet}$ and $\lambda_j\in\Lambda\cap Q^{a_i}$, so $i\le j$ by ~\eqref{eq:lambda-basis}.   We want to show that $i$ is the largest active index in $a_{\bullet}$ with $i\le j$.  Indeed, if $l$ is active in $a_{\bullet}$ with $i<l\le j$, then $\lambda_j\in P^{a_l}\cap Q^{b_{i'}}$.    But now $a_l > a_i$, so $a_l + b_{i'} > a_i + b_{i'} =  t+t'=d$.  Therefore $P^{a_l}\cap Q^{b_{i'}}=0$, contradiction. A similar argument shows $n(j) = t'$.  Therefore,
$$\codim_H (P^{m(j)} + Q^{n(j)} + \Lambda) = \codim_H (P^{t} + Q^{t'} + \Lambda)=1,$$
since $P^{t} + Q^{t'} $  is a hyperplane in $H$, and $\Lambda$ is contained in it by assumption.

It remains to show that if the conditions in the statement of Corollary~\ref{cor:at} do not all hold, then $\dim T_\Lambda  (\Sigma_{P^{\bullet},a_{\bullet}} \cap \Sigma_{Q^{\bullet},b_{\bullet}})= \rho-1.$  We prove the contrapositive. Suppose that $\dim T_\Lambda  (\Sigma_{P^{\bullet},a_{\bullet}} \cap \Sigma_{Q^{\bullet},b_{\bullet}})>\rho-1.$  By Theorem~\ref{thm:workhorse}, there is an index $j$ such that $\codim_H (P^{m(j)} + Q^{n(j)} + \Lambda) >0$.  Again, given that $\lambda_j \in P^{m(j)} \cap Q^{n(j)}$ and that $P^\bullet$ and $Q^\bullet$ are almost-transverse, it follows that either $m(j)+n(j)<d$ or that $m(j)=t$ and $n(j)=t'$.  But $m(j)+n(j)<d$ would imply $P^{m(j)} + Q^{n(j)} = H$, contradicting the codimension statement. So $m(j)=t$ and $n(j)=t'$, implying that $t=a_i$ and $t'=b_{i'}$ for active indices $i$ and $i'$ in $a_{\bullet}$ and $b_{\bullet}$ respectively.
(It is not possible that $m(j) = 0$ or $n(j) = 0$, since $\codim_H P^{m(j)}+Q^{n(j)}+\Lambda>0$.)
Furthermore, 
$$\langle \lambda_j \rangle = P^t\cap Q^{t'} \subseteq \Lambda \subseteq P^t + Q^{t'}$$
where the last containment holds again by the codimension assumption. 

Summarizing, we have shown that the {\em only} way that $\dim T_\Lambda  (\Sigma_{P^{\bullet},a_{\bullet}} \cap \Sigma_{Q^{\bullet},b_{\bullet}})> \rho-1$ is for all the conditions in the statement of Corollary~\ref{cor:at} to hold, in which case we have already proved that the dimension is exactly $\rho$.
\end{proof}


\begin{rem}\label{rem:arbitrary-flags}

For arbitrary flags $P^{\bullet}$ and $Q^{\bullet}$ and $\Lambda\in\Sigma^\circ_{P^\bullet,a_{\bullet}}\cap \Sigma^\circ_{Q^\bullet,b_{\bullet}}$, let $\tau\in S_d$ be 
the associated permutation from Proposition 
\ref{prop:sigma} (maintaining other notation as in Theorem 
\ref{thm:workhorse}). Then the extent to which 
the dimension of the tangent space at $\Lambda$ of 
$\Sigma_{P^{\bullet},a_{\bullet}} \cap \Sigma_{Q^{\bullet},b_{\bullet}}$
exceeds $\rho-1$ can be bounded in terms of $\tau$ as follows. We have:
\begin{equation} \label{eq:coxeter}
\dim T_\Lambda (\Sigma_{P^{\bullet},a_{\bullet}} \cap \Sigma_{Q^{\bullet},b_{\bullet}})
 \leq (\rho-1) + \inv( \omega \tau ).
\end{equation}
Here $\omega$ denotes the decreasing permutation $(d,d-1,\cdots,1)$, and $\inv(\omega\tau)$
denotes the inversion number of $\omega\tau$, i.e.\ the number of $i<j$ with 
$\omega\tau(i)>\omega\tau(j)$.\footnote{If one views $S_d$ as a Coxeter group
with reflections being the adjacent transpositions, then $\inv(\omega \tau)$
is also the Coxeter length of $\omega \tau$.}
We briefly sketch a proof of this more general inequality.

Using the second part of Proposition \ref{prop:sigma}, it follows that for 
each $i,j$, we have
$$ \dim P^i \cap Q^j = \# \{ i' \geq i:\ \tau(i') \geq j \}. $$
From this it follows that $\dim P^i \cap Q^j > \dim P^{i+1} \cap Q^j$ if and 
only if $\tau(i) \geq j$. Now, since we have 
$P^{a_i} \cap Q^{b_{\sigma(i)}} \neq P^{a_i+1} \cap
Q^{b_{\sigma(i)}}$, we find 
that $b_{\sigma(i)} \leq \tau(a_i)$. Note also that for all $j$, $m(j) \leq a_j$ and $n(j) \leq b_{\sigma(j)}$. Then:

\begin{eqnarray*}
\sum_{j=0}^r \codim ( \Lambda + P^{m(j)} + Q^{n(j)} ) &\leq&
\sum_{j=0}^r \codim(P^{m(j)} + Q^{n(j)})\\
&\leq& \sum_{j=0}^r \codim(P^{a_j} + Q^{b_{\sigma(j)}})\\
&\leq& \sum_{j=0}^r \codim( P^{a_j} + Q^{\tau(a_j)} )\\
&\leq& \sum_{j=0}^{d-1} \codim( P^j + Q^{\tau(j)}).
\end{eqnarray*}

Using that
$\dim(P^j \cap Q^{\tau(j)})=\#\{j' \geq j: \tau(j') \geq \tau(j)\}$, we 
compute that $\displaystyle \sum_{j=0}^{d-1} \codim( P^j
+ Q^{\tau(j)}) = \inv(\omega \tau)$, and the inequality \eqref{eq:coxeter} 
follows from \eqref{eq:the-answer}. 
When $P^{\bullet}$ and $Q^{\bullet}$ are almost-transverse, $\inv(\omega
\tau) = 1$, and 
the precise statement in Corollary
\ref{cor:at} can be deduced from characterizing the equality cases of the four inequalities above in the case where $\omega \tau$ is equal to an
adjacent transposition.
\end{rem}

\section{Linear series in positive genera}\label{sec:proofs}

We begin with a proposition that will show that our smoothness result, Theorem~\ref{thm:main} to be proved below, is sharp.

\begin{prop}\label{prop:non-smooth}
In the situation of Theorem \ref{thm:bg}, 
every point of
$G^r_d(X,(P,a_{\bullet}),(Q,b_{\bullet}))$ in the complement of
$G^{r,\circ}_d(X,(P,a_{\bullet}),(Q,b_{\bullet}))$ is singular.
\end{prop}

\begin{proof} This is a consequence of the standard construction of
the space $G^r_d(X,(P,a_{\bullet}),(Q,b_{\bullet}))$: we let 
$\wt{\sL}$ be a Poincar\'e line bundle on $X \times \Pic^d(X)$.
Take a sufficiently ample\footnote{Precisely, of degree strictly greater
than $2g-2$} effective divisor $D$ on $X$ with support
disjoint from $P$ and $Q$, and write $\wt{D}=D \times \Pic^d(X)$. Writing $p\colon X\times \Pic^d(X)\to\Pic^d(X)$ for projection, let
$G$ be the relative Grassmannian 
$\Gr(r+1,p_{*} \wt{\sL}(\wt{D}))$, equipped with structure map $\pi\colon G\to \Pic^d(X)$.  Let $\wt{\sV} \hookrightarrow \pi^* p_{*} \wt{\sL}(\wt{D})$ denote the universal subbundle.
Then $G^r_d(X)$ is cut out in $G$ by the condition that the induced map
$$\wt{\sV} \hookrightarrow \pi^*p_{*}\left( \wt{\sL}(\wt{D})|_{\wt{D}}\right)$$
vanishes identically. Because we have chosen $D$ to have support disjoint  
from $P$ and to be sufficiently ample, the space 
$G^r_d(X,(P,a_{\bullet}))$ is cut out by imposing the
additional Schubert condition that the maps
$$\wt{\sV} \hookrightarrow \pi^* p_{*}\left( \wt{\sL}(\wt{D})|_{a_j P}\right)$$
have rank at most $j$ for each $j$. 
Imposing the analogous condition at 
$Q$, we obtain $G^r_d(X,(P,a_{\bullet}),(Q,b_{\bullet}))$ as an 
intersection of three conditions: a determinantal condition (in fact a
complete intersection), and two relative Schubert cycles. It is routine
to check that for
$G^r_d(X,(P,a_{\bullet}),(Q,b_{\bullet}))$ to have dimension
$\rho(g,r,d,a_{\bullet},b_{\bullet})$, as asserted by Theorem \ref{thm:bg}, 
these three
conditions must intersect in the maximal codimension. Given that we know
from Theorem \ref{thm:bg} that
$G^r_d(X,(P,a_{\bullet}),(Q,b_{\bullet}))$ does in fact have dimension
$\rho(g,r,d,a_{\bullet},b_{\bullet})$, it then further
follows that in order for
$G^r_d(X,(P,a_{\bullet}),(Q,b_{\bullet}))$ to be smooth at any point,
that point must lie in the smooth locus of each of the three conditions,
and in particular of the two Schubert cycles. 

But we claim that
$G^{r,\circ}_d(X,(P,a_{\bullet}),(Q,b_{\bullet}))$ consists precisely of the
points of $G^r_d(X,(P,a_{\bullet}),(Q,b_{\bullet}))$ which lie in the 
smooth locus of both relative Schubert cycles. 
Indeed, each relative Schubert cycle is nothing but a locally constant family of Schubert varieties over the base $G$, so we are done by the standard characterization of the smooth locus of a Schubert variety (see the second part of Proposition \ref{prop:tangent-schubert}). 
\end{proof}

We now use our calculations in Grassmannians in \S2 to complete the proof of our 
main theorem, beginning with the case of genus $1$ in 
Theorem~\ref{thm:genus-1}.

\begin{proof}[Proof of Theorem \ref{thm:genus-1}] 
Set $\rho=\rho(1,r,d,a_{\bullet},b_{\bullet})$.
We may assume $d>0$, as otherwise the result is trivial. Thus,
$G^r_d(X)$ is a Grassmannian bundle over $\Pic^d(X)$, with the fiber
over a line bundle $\sL$ being canonically identified with 
$\Gr(r+1,\Gamma(X,\sL))\cong \Gr(r+1,d)$. The condition imposed by
requiring vanishing sequence at least $a_{\bullet}$ at $P$ then gives
a Schubert cycle in each fiber, corresponding to the complete flag determined by
vanishing order at $P$. The codimension of spaces in the flag corresponds
precisely to vanishing order except over the point $\sL \cong \sO_X(dP)$,
where no sections vanish to order precisely $d-1$, and vanishing to order
$d$ imposes codimension only $d-1$. Consequently, there are
two possibilities for $G^r_d(X,(P,a_{\bullet}))$. First, if $a_r<d$, it is
a relative Schubert cycle of codimension $\sum_j (a_j-j)$ in
$G^r_d(X)$, Cohen-Macaulay and flat over $\Pic^d(X)$.
Or,
if $a_r=d$, it is supported entirely over $\sL \cong \sO_X(dP)$ (even
scheme-theoretically), and is still a Schubert cycle, but of codimension
$(\sum_j (a_j-j))-1$. The same analysis applies to 
$G^r_d(X,(Q,b_{\bullet}))$, so we find that every fiber of 
\begin{equation}\label{eq:str-map} 
G^r_d(X,(P,a_{\bullet}),(Q,b_{\bullet})) \to \Pic^d(X)
\end{equation} 
is an intersection of a pair of Schubert cycles. The basic properties of
the map \eqref{eq:str-map} are analyzed for instance in Lemma 2.1 of
\cite{os18} and Proposition 2.1 of \cite{os26}; we review the main points
of this analysis in order to carry out the necessary tangent space analysis.

First, we see that in most fibers of \eqref{eq:str-map}, the relevant 
Schubert cycles are associated to transverse flags:
the only way in which the flags fail to be transverse is if 
$\sL \cong \sO_X(aP+(d-a)Q)$ for some $a\in\{1,\ldots,d-1\}$, which is 
unique by genericity of $P$ and $Q$; then the conditions of
vanishing to order $a$ at $P$ and $d-a$ at $Q$ intersect in dimension $1$
instead of dimension $0$. Thus, on fibers of \eqref{eq:str-map}
over points not of the form 
$\sO_X(aP+(d-a)Q)$ for $0 \leq a \leq d$, we have that the Schubert indexing
matches vanishing sequences, and the flags are transverse, so
the standard theory (see for instance Corollary \ref{cor:transverse})
gives us that (on these fibers) the space
$G^{r,\circ}_d(X,(P,a_{\bullet}),(Q,b_{\bullet}))$ is smooth (of 
relative dimension $(r+1)(d-r-1)-\sum_j (a_j-j)-\sum_j (b_j-j) = \rho-1$) over
$\Pic^d(X)$, and hence smooth of relative dimension $\rho$ 
over $\Spec k$. Similarly, it is easily
verified that we still obtain Richardson varieties over $\sO_X(aP+(d-a)Q)$ 
for $0<a<d$
unless $a$ occurs in $a_{\bullet}$ and $d-a$ occurs in $b_{\bullet}$.
Thus, we obtain the desired statement in these cases.
On the other hand, if $a=0$ or $a=d$, we have transverse intersection of
flags, but a potential difference in indexing. In the case $a=d$, the
difference in indexing arises only in that imposing vanishing order $d$ at 
$P$ is a codimension $d-1$ condition. Thus, this can only affect the final
term in the vanishing sequence, which is irrelevant for determining 
membership in 
$G^{r,\circ}_d(X,(P,a_{\bullet}),(Q,b_{\bullet}))$. We conclude 
that---whether or not $a_r=d$---the fiber of 
$G^{r,\circ}_d(X,(P,a_{\bullet}),(Q,b_{\bullet}))$ over
$\sO_X(dP)$ precisely corresponds to the open subset addressed
in Corollary \ref{cor:transverse}, and hence has the desired smoothness
property. The case that $a=0$ is the same, with $Q$ in place of $P$.

It thus remains to analyze the fibers with $\sL=\sO_X(aP+(d-a)Q)$, for 
$0<a<d$, and with $a$ 
occurring in $a_{\bullet}$ and $d-a$ occurring in $b_{\bullet}$. 
Our hypothesis on $P-Q$ implies that 
$\sO_X(aP+(d-a)Q) \not \cong \sO_X(a'P+(d-a')Q)$ for any $a \neq a'$,
so even in these cases, our flags in $\Gamma(X,\sL)$ are almost-transverse.
In the case that $a=a_j$ and $d-a=b_{r-j}$ for some $j$, then
$G^r_d(X,(P,a_{\bullet}),(Q,b_{\bullet}))$ is supported 
(scheme-theoretically) in the given fiber, and one checks (see 
the proof of Proposition 2.1 of \cite{os26}) that the nonempty fiber can still 
be described as a Richardson variety, by replacing $a_j$ with $a-1$, and
changing the choice of codimension-$a$ subspace in the first flag. Because
of this modification to the flag, we will still have that
$G^{r,\circ}_d(X,(P,a_{\bullet}),(Q,b_{\bullet}))$ corresponds precisely to 
the open subset treated in Corollary \ref{cor:transverse}, and is hence 
smooth (but this time of dimension $\rho$).

Finally, we consider the case that $a_j+b_{r-j}<d$ for all $j$, but we
have $\sL=\sO_X(aP+(d-a)Q)$, and $a=a_j$ and $d-a=b_{j'}$ for some $j,j'$ 
with $j+j'>r$; in particular, we must have $0<a<d$. See Example~\ref{ex:0202} for what is essentially the smallest nontrivial example of this case.
In this situation, 
the given fiber of
$G^{r,\circ}_d(X,(P,a_{\bullet}),(Q,b_{\bullet}))$ over $\Pic^d(X)$ may
be singular or even reducible, but at least it is pure of dimension
$\rho-1$: see the proof of Proposition 2.1 of \cite{os26}.
Moreover, as we have observed, the fiber is an intersection
of Schubert cycles associated to almost-transverse flags, so we can
invoke Corollary \ref{cor:at} to conclude that the fiber of
$G^{r,\circ}_d(X,(P,a_{\bullet}),(Q,b_{\bullet}))$ has tangent space
dimension equal to $\rho-1$ or $\rho$ everywhere. Moreover, the latter
occurs precisely at linear series $(\sL,V)$ satisfying the following
conditions: 
\begin{enumerate}
\item $V$ contains
a section $s$ vanishing to order $a$ at $P$ and $d-a$ at $Q$; 
\item $V$ is
contained in the linear span of the spaces of sections vanishing to order $a$ at
$P$ and order $d-a$ at $Q$; 
\item there is some $j$ with $a=a_j$ and $j$ active
in $a_{\bullet}$ in the sense of Definition \ref{defn:active}; 
\item there is some 
$j'$ with $d-a=b_{j'}$ and $j'$ active in $b_{\bullet}$.
\end{enumerate}

Thus, in order to complete the proof of the
theorem, we will prove that the space
$G^{r,\circ}_d(X,(P,a_{\bullet}),(Q,b_{\bullet}))$ is smooth of dimension
$\rho$ at every point of the given fiber by showing that
if $(\sL,V)$ is a point at which the tangent space of the fiber has
dimension $\rho$, then every tangent vector of the total space at $(\sL,V)$
is in fact vertical.
Accordingly, given $(\sL,V)$ satisfying the four conditions above, suppose $(\wt{\sL},\wt{V})$ is a first-order deformation of
$(\sL,V)$, and let $s \in V$ be a section vanishing to order $a$ at $P$ and
$d-a$ at $Q$.  We claim that $s$ has a lift $\tilde{s} \in \wt{V}$ which vanishes 
(scheme-theoretically) to order $a$ at $P$.  
Indeed, in the notation of the proof of Proposition~\ref{prop:non-smooth}, recall that on $G^{r}_d(X,(P,a_{\bullet}),(Q,b_{\bullet}))$ we have a map of vector bundles 
\begin{equation}\label{eq:rank-at-most-j}
\Phi\colon\tilde{\sV}\hookrightarrow\pi^*p_*\left(\tilde{\sL}(\tilde{D})|_{a_jP}\right)
\end{equation}
which has rank at most $j$.  The assumption that $(\sL,V)$ lies 
in the open subset $G^{r,\circ}_d(X,(P,a_{\bullet}),(Q,b_{\bullet}))$, 
together with the 
fact that $j$ is active in $a_{\bullet}$, says that on every point in an open neighborhood of $(\sL,V)$, the map~\eqref{eq:rank-at-most-j} has rank {\em exactly} $j$. In this situation, a standard argument shows that $\operatorname{ker}\Phi$ is locally free, and that the restriction map $(\operatorname{ker}\Phi)|_{(\tilde{\sL},\tilde{V})} \rightarrow (\operatorname{ker}\Phi)|_{(\sL,V)}$ is surjective. A sketch of this standard argument is as follows. The cokernel of $\Phi$ must be locally free (see e.g.~\cite[\S16.7]{ei1}), and hence also the image and kernel. Therefore, the short exact sequences $0\to\operatorname{ker}\Phi\to \tilde{\sV} \to \operatorname{im}\Phi\to0$ and $0\to\mathrm{im}\Phi\to\pi^*p_*\left(\tilde{\sL}(\tilde{D})|_{a_jP}\right) \to\operatorname{cok}\Phi\to 0$ remain exact after base change, and hence taking kernels commutes with base change. 

From surjectivity of the restriction map above, we deduce that our section $s$ admits a lift $\tilde{s}$ that also vanishes (scheme-theoretically) to order $a$ at $P$.
Similarly, $s$ must have
another lift which vanishes 
(scheme-theoretically) to order $d-a$ at $Q$; since it is another lift of $s$, it can be 
expressed as $\tilde{s}+\epsilon v$ for some $v \in V$.

Now, recall our hypothesis that $V$ is contained in the 
span of sections vanishing to order at least $a$ at $P$ and at least $d-a$ at $Q$.
Write $v=v_1 + v_2$, where $v_1$ vanishes to order at least $a$ at 
$P$ and $v_2$ vanishes to order at least $d-a$ at $Q$. But then 
$\tilde{s}+ \epsilon v_1$ still vanishes to order $d-a$ at $Q$, and also vanishes to 
order $a$ at $P$. This forces $\wt{\sL}$ to be the trivial deformation of $\sL$,
yielding the desired verticality
assertion and the theorem.
\end{proof}

%
To conclude the proof of our main theorem, we need to make use of the
Eisenbud-Harris theory of limit linear series. We first set up notation
for our reducible curves, and recall the relevant definitions.

\begin{sit}\label{sit:chain} Fix $g,d,n$. Let $Z_1,\dots,Z_n$
be smooth projective curves, with (distinct) points $P_i, Q_i$ on $Z_i$ for each
$i$, and let $X_0$ be the nodal curve obtained by gluing $Q_i$ to $P_{i+1}$ 
for $i=1,\dots,n-1$.
\end{sit}

\begin{defn}\label{defn:lls} Given $r,d$, a \textbf{limit linear series}
of dimension $r$ and degree $d$ on $X_0$ consists of a tuple $(\sL^i,V^i)$ of 
linear series of dimension $r$ and degree $d$ on the $Z_i$, satisfying the
following condition: if $a^i_{\bullet},b^i_{\bullet}$ are the vanishing
sequences of $(\sL^i,V^i)$ at $P_i$ and $Q_i$ respectively, then we require
\begin{equation}\label{eq:eh-ineq} 
b^i_j+a^{i+1}_{r-j} \geq d
\end{equation}
for all $i=1,\dots,n-1$ and $j=0,\dots,r$. If \eqref{eq:eh-ineq} is an
equality for all $i,j$, we say that the limit linear series is
\textbf{refined}.

The space of all such limit linear series on $X_0$ is denoted by
$G^r_d(X_0)$. If we have sequences $a_{\bullet}$ and $b_{\bullet}$,
we also have the closed subscheme 
$G^r_d(X_0,(P_1,a_{\bullet}),(Q_n,b_{\bullet})) \subseteq G^r_d(X_0)$
consisting of limit linear series such that, following the above notation,
we have $a^1_{\bullet} \geq a_{\bullet}$ and $b^n_{\bullet} \geq b_{\bullet}$.
Finally, denote by
$G^{r,\circ}_d(X_0,(P_1,a_{\bullet}),(Q_n,b_{\bullet})) \subseteq 
G^r_d(X_0,(P_1,a_{\bullet}),(Q_n,b_{\bullet}))$ the open subscheme
consisting of refined limit linear series which further satisfy
\begin{eqnarray*}
\dim V^1(-a_j P_1) =r+1-j &\text{for $j>0$ active in $a$,}\\
\dim V^n(-b_j Q_n) =r+1-j &\text{for $j>0$ active in $b$.}
\end{eqnarray*}
We comment that these last two conditions can be re-expressed purely in terms of $a^1_\bullet$ and $b^n_\bullet$ as follows: for all $j>0$ active in $a$, we require
$\#\{j':a^1_{j'} \ge a_j\} = r+1-j$, and similarly for $b$.
\end{defn}


We are now ready to prove our main smoothness result. 

\begin{proof}[Proof of Theorem \ref{thm:main}]
In Situation \ref{sit:chain}, observe that we can decompose the limit
linear series space
$G^{r,\circ}_d(X_0,(P_1,a_{\bullet}),(Q_n,b_{\bullet}))$ into disjoint open
subsets according to the vanishing sequences at each node, i.e.~according
to the possible values in the left hand side of the equalities~\eqref{eq:eh-ineq}.
Then each such open subset is almost a product over $i$ of spaces of the form
$G^{r,\circ}_d(Z_i,(P_i,a^i_{\bullet}),(Q_i,b^i_{\bullet}))$. In fact, it
is an open subset of this product, since the refinedness condition completely
fixes the vanishing sequences at the nodes. If further
each $Z_i$ has genus $0$ or $1$, and for the $Z_i$ of genus $1$ we
suppose that $P_i-Q_i$ is not $m$-torsion for any $m \leq d$, then we know 
that each $G^{r,\circ}_d(Z_i,(P_i,a^i_{\bullet}),(Q_i,b^i_{\bullet}))$ is
smooth. Indeed, the genus-$1$ case is Theorem \ref{thm:genus-1}, while the
genus-$0$ case is well known, but follows in particular immediately from
Corollary \ref{cor:transverse} taking into account that 
$\dim \Gamma(\PP^1,\sO(d))=d+1$, so there is a shift of $1$ in the value
of $d$. We thus conclude that
$G^{r,\circ}_d(X_0,(P_1,a_{\bullet}),(Q_n,b_{\bullet}))$ is also smooth,
of dimension $\rho(g,r,d,a_{\bullet},b_{\bullet})$.

Now, fix $n=g$ and suppose each $Z_i$ has genus $1$.
Let $B$ be the spectrum of a discrete valuation ring, and $\pi:X \to B$ be
a flat, proper family of curves of genus $g$, with 
$X$ regular, the generic fiber $X_{\eta}$ smooth, and the special fiber
isomorphic to $X_0$. Further assume that $\pi$ has sections $P$, $Q$,
specializing to $P_1$ and $Q_n$ respectively on $X_0$.

Suppose that we have a closed point of 
$G^{r,\circ}_d(X_{\eta},(P_{\eta},a_{\bullet}),(Q_{\eta},b_{\bullet}))$.
After a base change, we may assume that the corresponding linear series is defined
on $X_{\eta}$, and that all ramification points are also
rational over the base field. Blow up the nodes in $X_0$ as necessary
to resolve any singularities of $X$ resulting from the base change,\footnote{The preceding constitutes
an alternative for the argument of Theorem 2.6 of \cite{e-h1}, avoiding 
invocation of the stable reduction theorem.} and finally, blow up $P_1$ and
$Q_n$ as necessary so that no generic ramification point distinct from
$P$ or $Q$ limits to $P_1$ or $Q_n$ in the special fiber.
Denote the
resulting family by $\pi':X' \to B'$, and the special fiber by $X'_0$,
and write $P'$ and $Q'$ (respectively, $P'_1$ and $Q'_n$) for the
resulting sections of $\pi'$ and their restrictions to $X'_0$. 
Then $X'_0$ is obtained by $X_0$ by base extension and insertion of 
chains of genus-$0$ curves at the nodes and at $P_1$ and $Q_n$.
By construction, none of the ramification points on $X_\eta$ can specialize to nodes of 
$X'_0$, so by Proposition 2.5 of \cite{e-h1} (and using the characteristic
$0$ hypothesis), the extension of the given linear series is a refined limit
linear series on $X'_0$. Moreover, by the same argument, the ramification
at $P'_1$ and at $Q'_n$ must be precisely equal to the ramification at
$P_{\eta}$ and $Q_{\eta}$, so
that the induced limit linear series lies in 
$G^{r,\circ}_d(X'_0,(P'_1,a_{\bullet}),(Q'_n,b_{\bullet}))$.
But as we have discussed above, this space is smooth. According to
Proposition 2.7 of \cite{c-h-t1} (see also Theorem 3.4 of \cite{os26} 
for the situation with imposed ramification) there is a flat relative moduli 
space recovering linear series on the generic fiber and (refined) limit 
linear series on the special fiber.\footnote{Proposition 2.7 of 
\cite{c-h-t1} claims incorrectly to apply to all limit linear series in
the special fiber, but it only applies to the refined limit linear series,
since it uses the construction of Eisenbud and Harris
\cite{e-h1}, which contains only the refined limit linear series. However,
this is enough for our purposes.}
Alternatively, the more recent construction of
\cite{o-m1} gives an independent approach which extends to all limit linear
series.
Either way, it follows that the original point of 
$G^{r,\circ}_d(X_{\eta},(P_{\eta},a_{\bullet}),(Q_{\eta},b_{\bullet}))$
must have been smooth as well.

Now, since the spaces we are considering are in general not proper, the condition that $G^{r,\circ}_d$ is smooth is not open in families. However, the condition does define a constructible subset of $\cM_{g,2}$, 
and the generic fibers of the possible families $\pi$ as above correspond to
a Zariski-dense subset, so we conclude the main smoothness statement of the 
theorem. The fact that the remaining points are not smooth is Proposition 
\ref{prop:non-smooth}. 

The statement on codimension of singularities follows 
from the observation that a point in the complement of 
$G^{r,\circ}_d(X,(P,a_{\bullet}),(Q,b_{\bullet}))$ 
is a union of closed subvarieties of the form $G^{r}_d(X,(P,a'_{\bullet}),(Q,b'_{\bullet}))$, where $a'_\bullet \ge a$ and $b'_\bullet \ge b$ are sequences such that for some $j>0$ active in $a_\bullet$, we have $a'_{j-1} \ge a_j$, implying that $a'_{j-1}\ge a_{j-1}+2$ and $a_j' \ge a_j+1$ and hence $\sum a'_\bullet \ge \sum a_\bullet+3$; or analogously for $b_\bullet$.
We then conclude normality from the Cohen-Macaulayness and
Serre's criterion, and the irreducibility statement follows immediately
from the connectedness in the case $\widehat{\rho} \geq 1$.
\end{proof}

\begin{rem}\label{rem:hook}
The singular locus of $G^r_d(X,(P,a_\bullet))$ can also be described as follows.  Given $j>0$ with $a_j > a_{j-1}+1$, let $m\ge 0$ be largest such that $a_{j+m} = a_j + m$. 
Let 
\begin{align*}
a^{(j)}_\bullet &= (a_0,\ldots,a_{j-2},a_{j},a_j\!+\!1,\ldots,a_{j+m}\!+\!1,\ldots,a_r).
\end{align*}
\noindent Then $$G^r_d(X,(P,a_\bullet))^\text{sing} = \!\!\bigcup_{\substack{j>0\\a_j>a_{j-1}+1}} \!\!G^r_d(X,(P,a^{(j)}_\bullet)).$$
In terms of Young diagrams, let $\alpha_\bullet = (a_i\!-\! i)_{i=0}^r$ and let $\operatorname{Young}(a_\bullet)$ be the complement of the partition $\alpha_\bullet$ in a $(r\!+\!1)\times (d\!+\!r)$ box.  Then $\{\operatorname{Young}(a^{(j)}_\bullet)\}$ consists of all diagrams obtained from $\operatorname{Young}(a_\bullet)$ by removing hooks.  

In the twice-marked case, it is immediate from Theorem~\ref{thm:main} that the singular locus of $G^r_d(X,(P,a_\bullet), (Q,b_\bullet))$ is
\begin{align*}
\left(G^r_d(X,(P,a_\bullet))^\text{sing} \cap G^r_d(X, (Q,b_\bullet))\right)  
&\cup  
\left(G^r_d(X,(P,a_\bullet)) \cap G^r_d(X, (Q,b_\bullet))^\text{sing}\right).
\end{align*}
These descriptions generalize the well-known characterization of Grassmannian Schubert singularities, and Richardson singularities, in terms of hook removals \cite[Theorem 9.3.1]{b-l2}.
\end{rem}

\begin{ex}\label{ex:0202}
We provide here essentially the smallest interesting example in the $g=1$ case, exhibiting a fiber of $G^r_d(X, (P, a_\bullet), (Q,b_\bullet))\to \Pic^d(X)$ that is not a Richardson variety, and is in fact reducible.  Let $a_\bullet = b_\bullet = (0,2)$, let $\sL = \sO_X(2P + 2Q)$, and consider the fiber of $G^1_4(X, (P, a_\bullet), (Q,b_\bullet))$ over $\sL$.  This fiber has two irreducible components $Z_1$ and $Z_2$, each isomorphic to $\mathbb{P}^2$, meeting along a $\mathbb{P}^1$. It may be described as the variety of lines in $\mathbb{P}^3$ that meet two fixed lines that themselves intersect at a point.  

In this situation, the vertical tangent space at a point in $Z_1\cap Z_2\cong \mathbb{P}^1$ has dimension jumping up to $3$.  Now, the fiber of $G^{1,\circ}_4(X, (P, a_\bullet), (Q,b_\bullet))$ over $\sL$ is obtained from that of $G^1_4(X, (P, a_\bullet), (Q,b_\bullet))$ by removing two points of $Z_1\cap Z_2$. Those two points correspond to the space of sections of $\sL$ vanishing to order at least $2$ at $P$, respectively the space of sections of $\sL$ vanishing to order at least $2$ at $Q$.  
Then Theorem~\ref{thm:genus-1} asserts that on $Z_1\cap Z_2$, except for at those two points, $G^1_4(X, (P, a_\bullet), (Q,b_\bullet))$ has no horizontal tangent vectors.
\end{ex}

\begin{rem}\label{rem:wrong}
Proposition 3.2 of \cite{c-h-t1} asserts that for a general marked curve,
the space $G^r_d(X,(P,a_{\bullet}))$
is connected when $\rho>0$, and is smooth at all points corresponding to
complete linear series. However, Proposition \ref{prop:non-smooth} as well
as the results in \cite{os26} refute both of these statements; see
Example \ref{ex:wrong} below for counterexamples. No argument is 
provided in \cite{c-h-t1} for connectedness, but an argument is provided
for the smoothness assertion. The authors have indicated to us that
upon further examination, their argument appears correct but implicitly uses 
that $h^0(L(-a_i P))=r+1-i$ for all $i$, which is consistent with the 
criterion we provide in Theorem \ref{thm:main}.
\end{rem}

\begin{ex}\label{ex:wrong}
The space $G^r_d(X,(P,a_{\bullet}))$ need not be smooth on the locus of
complete linear series. For instance,
suppose that $g \geq 11$ is odd, and set
$d=\frac{g+1}{2}+3$, and $r=1$, with $a_{\bullet}=0,2$. Then $\rho=4$,
and we see that $G^r_d(X,(P,a_{\bullet}))$ consists entirely of complete
linear series when $X$ is general, since $\rho$ becomes negative if
$r$ is increased to $2$. However, the locus of linear series with a 
double base point at $P$ is nonempty, and is contained in the singular
locus by Proposition \ref{prop:non-smooth}.

The assertion that $G^r_d(X,(P,a_{\bullet}))$ is connected
is true when $\widehat{\rho}>0$ by the main 
theorem of \cite{os26}, but it can fail when $\widehat{\rho}=0$ even with
a single ramification point. For instance, consider the case that $g=3$,
$r=2$, $d=6$, and $a_{\bullet}=0,3,5$. According to Remark 4.6 of
\cite{os26}, we can compute the number of connected components of
$G^r_d(X,(P,a_{\bullet}))$ for $X$ general by counting the number of
connected components of the space of limit linear series on our chain
of genus-$1$ curves. We see that in order to get a valid sequence 
(i.e., nonnegative and without repetitions) for $b^3_{\bullet}$, we 
must have $b^3_{\bullet}$ of the form $0,1,*$, and for
each $i$ we must have some $j$ with $a^i_j+b^i_{r-j}=d$, 
with $j=1$ occurring exactly once and $j=2$ occurring 
exactly twice. This leads to two possibilities: we could have
$a^2_{\bullet}=*,3,6$ and $a^3_{\bullet}=*,4,6$, or we could have
$a^2_{\bullet}=*,4,5$ and $a^3_{\bullet}=*,4,6$. We see that both
possibilities can occur in valid limit linear series, and that we obtain
two distinct connected components of the limit linear series space. Indeed,
in each case the underlying line bundles on each $Z_i$ are uniquely
determined, but they are distinct on $Z_1$ and $Z_2$. In the
first case we must have 
$\sO_{Z_1}(3P_1+3Q_1)$ on $Z_1$ and $\sO_{Z_2}(6P_2)$ on $Z_2$, while in
the second case we must have
$\sO_{Z_1}(5P_1+Q_1)$ on $Z_1$ and $\sO_{Z_2}(4P_2+2Q_2)$ on $Z_2$.
\end{ex}

\newcommand{\noopsort}[1]{} \newcommand{\printfirst}[2]{#1}
  \newcommand{\singleletter}[1]{#1} \newcommand{\switchargs}[2]{#2#1}
\providecommand{\bysame}{\leavevmode\hbox to3em{\hrulefill}\thinspace}
\providecommand{\MR}{\relax\ifhmode\unskip\space\fi MR }
\providecommand{\MRhref}[2]{%
  \href{http://www.ams.org/mathscinet-getitem?mr=#1}{#2}
}
\providecommand{\href}[2]{#2}


\end{document}